\documentclass{amsart}
\usepackage{amssymb}
\newtheorem{thm}{Theorem}[section]
\newtheorem{proposition}[thm]{Proposition}
\newtheorem{corollary}[thm]{Corollary}
\newtheorem{lemma}[thm]{Lemma}
\newcommand{\dontprint}[1]{\relax}
\newcommand{\bx}[1]
{\begin{picture}(8,8)\put(0.5,-1.0){\framebox(7,7){$\scriptstyle{#1}$}}\end{picture}}
\title[Polynomial solutions of the KZ equations]
{Polynomial solutions of the Knizhnik--Zamolodchikov equations and Schur--Weyl duality}
\author{Giovanni Felder}
\address{Department of mathematics, ETH Zurich,
8092 Zurich, Switzerland}
\email{giovanni.felder@math.ethz.ch}
\author{Alexander P. Veselov}
\address{Department of Mathematical Sciences,
Loughborough University, Loughborough LE11 3TU, UK  and Landau
Institute for Theoretical Physics, Moscow, Russia}
\email{A.P.Veselov@lboro.ac.uk}
\begin{document}
\maketitle

\begin{abstract}
An integral formula for the solutions of Knizhnik-Zamolodchikov
(KZ) equation with values in an arbitrary irreducible
representation of the symmetric group $S_N$ is presented for
integer values of the parameter. The corresponding integrals can
be computed effectively as certain iterated residues determined by
a given Young diagram and give polynomials with integer
coefficients. The derivation is based on Schur-Weyl duality and
the results of Matsuo on the original $SU(n)$ KZ equation. The
duality between the spaces of solutions with parameters $m$ and
$-m$ is discussed in relation with the intersection pairing in the
corresponding homology groups.
\end{abstract}

\section{Introduction}
Let $G$ be a finite Coxeter group, $R$ be the corresponding root
system,  $m_\alpha, \, \alpha \in R$ be a system of multiplicities,
which is a $G$-invariant function on $R.$ Let $W$ be an irreducible
representation of $G$ and define the {\em Knizhnik--Zamolodchikov
equation} (KZ equation) related to $W$ as the following system
\[
\partial_\xi \psi = \sum_{\alpha\in
R_+}m_\alpha\frac{(\alpha,\xi)}{(\alpha,z)}
(s_\alpha+1) \psi,
\]
where $s_{\alpha}$ are the corresponding reflections acting on
$W$-valued functions $\psi(z)$. If the multiplicities $m_\alpha$ are
integers then all the solutions of the corresponding systems are
homogeneous polynomials (see  \cite {O, FV}) of degree equal to the
value of the central element $\sum_{\alpha\in R_+} m_\alpha
(s_\alpha+1)$ in the irreducible representation $W$. The finding of
these solutions is an important part of the description of the
so-called $m$-harmonic polynomials \cite{FeiV}. In the paper
\cite{FV} these solutions were found explicitly in the simplest case
of the standard (reflection) representation of $G= S_N.$

The main result of the present paper is an explicit integral formula
for the solutions of the corresponding KZ equation
\begin{equation}
\label{KZ}
\partial_i \psi = m \sum_{j \neq i}^N \frac{s_{ij}+1}{z_i - z_j}
\psi, \quad i=1,\dots,N
\end{equation}
with values in an arbitrary irreducible representation of the
symmetric group $S_N$ for any positive integer $m.$ Our approach is
based on Schur-Weyl duality and the results of Matsuo, who found
some integral formulas for the solutions of the original $SU(n)$ KZ
equation \cite{M}.

The main construction is the following. Let $\lambda$ be the {\it
Young diagram} with $N$ boxes and $n$ rows of lengths
$\lambda_1,\dots,\lambda_n$ with $\lambda_1\geq\cdots\geq\lambda_n >
0$. It is well-known that for any such diagram one can construct an
irreducible representation $W^{\lambda}$ of the symmetric group
$S_N$  and any irreducible representation of $S_N$ can be described
in this way (see e.g. \cite{Fulton}).

The space $W^{\lambda}$ has a basis $v_T$ labeled by the set
$\mathcal T(\lambda)$ of all {\it standard tableaux} on $\lambda$,
which are the numberings $T\colon \lambda\to \{1,\dots,N\}$ of the
boxes of $\lambda,$ increasing from left to right and from top to
bottom.

The fundamental set of solutions of the KZ equation with values in
$W^\lambda$ can be also labeled by the set $\mathcal T(\lambda):$
\begin{equation}
\label{sol}
\psi_T (z_1, \dots, z_N) = \sum_{T' \in \mathcal T(\lambda)}
\psi_{T, T'} (z_1, \dots, z_N) v_{T'}.
\end{equation}
The components $\psi_{T, T'} (z_1, \dots, z_N)$ are known to be
polynomial in $z_1, \dots, z_N$ \cite{FV}. We can give now an
explicit formula for $\psi_{T, T'} (z_1, \dots, z_N)$ as an integral
\begin{equation}
\label{int}
\psi_{T, T'} = \int_{\sigma_T} \omega_{T'}
\end{equation}
of some rational differential form  $\omega_{T'}$ over a certain
cycle $\sigma_T$ in the top homology of the following configuration
space $C_{\lambda} (z_1, \dots, z_N)$ related to Young diagram
$\lambda.$

For given $\lambda = (\lambda_1, \dots, \lambda_n)$ let us define
the integers $m_i,\,  i=1, \dots, n$ from the relation $\lambda =
(m_0-m_1, m_1-m_2, \dots, m_{n-2}-m_{n-1}, m_{n-1}).$ Explicitly we
have
$$m_0= \lambda_1 + \lambda_2 + \dots + \lambda_n = N, \,\, m_1 =
\lambda_2 + \dots + \lambda_n, \,\, \dots, m_{n-2} = \lambda_{n-1} +
\lambda_n, m_{n-1} = \lambda_n,$$ so that $m_s$ is the number of
boxes in the rows of $\lambda$ strictly lower than $s.$

Consider  $n$ finite sets $X_0, X_1, \dots, X_{n-1}$ of points on
the complex plane $\mathbb{C}$ consisting of $m_0, \dots, m_{n-1}$
points respectively with the condition that $X_i$ and $X_{i+1}$ have
no common points for all $i = 1, \dots, n-2.$ Let us denote the
elements of $X_0$ as $z_1, \dots, z_N$ and fix them. The
corresponding configuration space of all admissible $\{X_1, \dots,
X_{n-1}\}$ is our $C_{\lambda} (z_1, \dots, z_N).$  It has the
dimension $$d_{\lambda} = \sum_{i=1}^{n-1} m_i = \sum_{r=1}^{n}
(r-1) \lambda_r,$$  and can be described as the following subset in
$\mathbb C^{d_{\lambda}}:$
\begin{equation}
\label{e-C-lambda}
C_{\lambda} (z_1, \dots, z_N) = \{ t^{b}_s \in
\mathbb C, b \in \lambda, 1 \le s \le r(b)-1: \,\, t^{b}_{s+1} \neq
t^{b'}_s, \, t^{b}_{1} \neq z_{k} \},
\end{equation} where we have denoted the elements of $X_s$  as $$X_s
= \{ t^b_s \in \mathbb C , b \in \lambda, \,  r(b)>s \},$$ and
$r(b)$ is the row which the corresponding box $b$ belongs to.

On this space $C_{\lambda} (z_1, \dots, z_N)$ we have a natural
action of the group $G_{\lambda} = S_{m_1} \times  S_{m_2} \times
\dots \times  S_{m_{n-1}}.$ The cycles $\sigma_T$ in the top
homology group $H_{\mathit{top}}(C_{\lambda} (z_1, \dots, z_N))$ can
be defined for any numbering $T$ of $\lambda$, namely any bijection
$T\colon \lambda\to \{1,\dots,N\}.$ Consider first the product
$\Gamma_T$ of circles consecutively surrounding anti-clockwise $z_k$
with the variables $t^{b}_1, \dots, t^{b}_{r(b)-1}, \, b =
T^{-1}(k)$ located on these circles:
\begin{equation}
\label{Gamma}
\Gamma_T = \{  t^b_s \in \mathbb C: | t^b_s - z_k | =
\epsilon s, \, b = T^{-1}(k)\}
\end{equation}
for any real positive $\epsilon$ small enough. The corresponding
cycle $\sigma_T$ is the skew-symmetrisation of $\Gamma_T$ by the
action of $G_{\lambda}$:
\begin{equation}
\label{sigma}
\sigma_T = \sum_{g \in G_{\lambda}} (-1)^{g}  g_*(\Gamma_T),
\end{equation}
where $(-1)^{g}$ denote the sign of $g,$ which is the product of
signs of the corresponding permutations in $S_{m_i}.$

The form $\omega_T$ has the form
\begin{equation}
\label{omega}
\omega_T = \frac {1}{(2 \pi i)^{d_{\lambda}}} \Phi_{\lambda}^m \phi_T dt,
\end{equation}
 where
\begin{equation}
\label{Phi}
\Phi_{\lambda} = \prod_{i<j}^N (z_i - z_j)^2
\prod_{s, b \neq b'} (t^{b}_s - t^{b'}_s)^2
\prod_{s, b, b'} (t^{b}_{s+1} - t^{b'}_s)^{-1}
\prod_{ k, b} (t^{b}_{1} - z_{k})^{-1} ,
\end{equation}
\begin{equation}
\label{phi}
\phi_T = \prod_{s, b} (t^{b}_{s+1} - t^{b}_s)^{-1}
\prod_{ b} (t^{b}_{1} - z_{T(b)})^{-1}
\end{equation}
and $dt = \prod_{s, b}  dt^{b}_s$ is the exterior product of the
differentials of all the coordinates (the order is not essential
since it is only changes sign).

\begin{thm}\label{t-main} For any given positive integer $m$,
the integral formulas (\ref{sol}), (\ref{int}) with the cycles
$\sigma_T$ and forms $\omega_T,  T \in \mathcal T(\lambda)$ defined
above give a basis in the space of solutions of the KZ equation
(\ref{KZ}) with values in the irreducible $S_N$-module
$W^{\lambda}.$ The integral (\ref{int}) can be effectively computed
as an iterated residue and gives a polynomial in $z_1, \dots, z_N$
with integer coefficients.
\end{thm}

We have derived these formulas from the results of Matsuo \cite{M}
using the Schur-Weyl duality. The fact that in such a way we get all
the solutions of the corresponding KZ equation does not follow from
\cite{M} and needs to be proven independently. Using representation
theory it turns out that it is sufficient to prove that one of the
integrals $\psi_{T,T'}$ does not vanish identically. This we prove
by evaluating $\psi_{T,T}$ for the tableau assigning $k$ to the
$k$th box, counted from left to right and from top to bottom, in the
asymptotic region $0\ll |z_1|\ll \cdots\ll |z_N|$. We find
\[
\psi_{T,T}(z)\sim C\, \prod_{b\in\lambda}
z_{T(b)}^{m(T(b)-1+c(b)-r(b))}+\cdots,
\]
for some integer $C\neq0$, see Lemma \ref{l-nonzero}, where
$c(b),r(b)$ are the coordinates (column and row number) of the box
$b\in\lambda$. As a by-product, we obtain a new derivation of the
interesting formula, due to Frobenius \cite{F}  (see \cite{Mac},
exercise 7 in Chapter I and comment on p.~134) for the value
$f_2(\lambda)$ of the central element $\sum_{i<j}s_{ij}$ in the
representation $W^\lambda$:
\[
f_2(\lambda)=\sum_{b\in\lambda}(c(b)-r(b)).
\]
Theorem \ref{t-main} applies to the case of positive integer $m$.
The case of negative $m$ can be reduced to it by using the
isomorphism between the space of solutions $\mathit{KZ}(V,m)$ of
the KZ equation with values in the representation $V$ and
parameter $m$ and the space
$\mathit{KZ}(V\otimes\mathrm{Alt},-m)$, where
$\mathrm{Alt}=\mathbb C\epsilon$ is the alternating
representation. Indeed it is easy to check that if
$\psi\in\mathit{KZ}(V,m)$ then $\phi=\prod_{i>j}(z_i-z_j)^{-2m}
\psi\otimes\epsilon\in\mathit{KZ}(V\otimes\mathrm{Alt},-m)$. In
particular it follows that for negative $m$ all solutions are
rational functions. In the last section we discuss also the
duality between $\mathit{KZ}(V,m)$ and $\mathit{KZ}(V^{*},-m)$
given by the canonical map $$\mathit{KZ}(V,m) \otimes
\mathit{KZ}(V^{*},-m) \rightarrow {\bf C}$$ in relation with the
intersection pairing in the corresponding homology groups.

The case we consider can be viewed as a very degenerate limit of
the general theory of hypergeometric solutions of KZ equations
associated with Kac--Moody algebras, see \cite{Varchenko} and
references therein. We should mention that similar integral
formulas were found also by Kazarnovski-Krol \cite{KK} in the
theory of zonal spherical functions of type $A_n,$ but
combinatorics of the corresponding configuration space is very
different and not related to Young diagrams.

\section{Schur--Weyl duality}

We start with the classical Schur-Weyl duality between the
representations of the general linear and symmetric groups.

Let $V$ be an $n$-dimensional complex vector space. Then the
symmetric group $S_N$ on $N$ letters acts on the $N$-fold tensor
product $V^{\otimes N}=V\otimes\cdots\otimes V$ by permutations of
factors and this action commutes with the diagonal action of
$GL(V)$. The Schur--Weyl theorem states that, as a $GL(V)\times S_N$
module, $V^{\otimes N}$ has a decomposition into a direct sum
\begin{equation}\label{e-Ribera}
V^{\otimes N}\cong\oplus_\lambda M^\lambda\otimes W^\lambda
\end{equation}
where $M^\lambda$ are inequivalent irreducible $GL(V)$-modules and
$W^\lambda$ are inequivalent irreducible $S_N$-modules. The sum is
over partitions of $N$ into at most $n$ parts, namely sequences of
integers $\lambda_1\geq\cdots\geq\lambda_n\geq 0$ with
$\sum\lambda_i=N$.

Moreover if $n\geq N$ all irreducible $S_N$ modules appear. Thus we
can realise every irreducible $S_N$-module as
\[
W^\lambda=\mathrm{Hom}_{GL(V)}(M^\lambda,V^{\otimes N}),
\]
for any $V$ of dimension $\geq N$.

An explicit description of $W^\lambda$ is obtained from the
description of $M^\lambda$ as a highest weight module and depends on
a choice of basis of $V$. Namely, let us fix a basis $e_1,\dots,e_n$
of $V$ and introduce the decomposition
$\mathfrak{gl}(V)=\mathfrak{n}_{-}\oplus \mathfrak{h}\oplus
\mathfrak{n}_+$ of the Lie algebra of $GL(V)$ into strictly lower
triangular, diagonal and strictly upper triangular $n\times n$
matrices. For any $\mathfrak{ gl}(V)$-module $E$ and any
$\mu\in\mathfrak h^*=\mathbb C^n$,  denote by $E_\mu=\{v\in
E\,|\,x\cdot v=\mu(x)v\}$ the weight subspace of weight $\mu$. The
space of {\em primitive vectors} of weight $\mu$ in $E$ is
\[
E_\mu^{\mathfrak n_+}=\{v\in M_\mu\,|\,a\cdot v=0, \forall
a\in\mathfrak n_+\}.
\]
The irreducible module $M^\lambda$ is uniquely characterised by
having a non-zero primitive vector $v_\lambda$ of weight $\lambda$,
which is unique up to normalization. Moreover $M^\lambda$ is
generated over $U(\mathfrak n_-)$ by $v_\lambda$ and thus
\begin{equation}\label{e-hwm}
M^\lambda=\mathbb C\,v_\lambda\oplus \mathfrak n_-\, M^\lambda.
\end{equation}
An element of $\mathrm{Hom}_{GL(V)}(M^\lambda,V^{\otimes N})$ is
then uniquely determined by the image of the primitive vector and we
obtain an isomorphism of $S_N$-modules
\[
W^\lambda= (V^{\otimes N})_\lambda^{\mathfrak n_+}.
\]
{}From this realization we obtain a basis of $W^\lambda$ labeled by
standard Young tableaux, making connection to the classical
construction of $W^\lambda$ as a Specht module, see \cite{Fulton}.
Here is how it works: let $\lambda$ also denote the Young diagram
with $N$ boxes with rows of lengths $\lambda_1,\dots,\lambda_m$. To
each numbering $T\colon \lambda\to \{1,\dots,N\}$ of the boxes we
associate a weight vector $e_T=e_{\alpha_1}\otimes\cdots\otimes
e_{\alpha_N}\in (V^{\otimes N})_\lambda$ so that $\alpha_k=i$
whenever $T^{-1}(k)$ is in the $i$th row. For example, if $T$ is the
numbering

\smallskip
\begin{picture}(100,60)(80,120)
\put(100,155){\framebox(14.50,14.60)[c]{3}}
\put(115,155){\framebox(14.60,14.60)[c]{4}}
\put(130,155){\framebox(14.60,14.60)[c]{1}}
\put(100,140){\framebox(14.60,14.60)[c]{2}}
\put(115,140){\framebox(14.60,14.60)[c]{5}}
\put(100,125){\framebox(14.60,14.60)[c]{6}}
\end{picture}
\smallskip

\noindent of $\lambda=(3,2,1)$, then $e_T=e_1\otimes e_2\otimes
e_1\otimes e_1\otimes e_2\otimes e_3$.

The symmetric group $S_N$ acts on the set of numberings of $\lambda$
and we have $\sigma e_T=e_{\sigma T}$ for any $\sigma\in S_N$. For
any numbering $T$ of $\lambda$, the row group $R(T)$ is the subgroup
of $S_N$ consisting of permutations preserving the image of each
row. Similarly, we have the column group $C(T)\subset S_N$. Two
numberings $T,T'$ give the same weight vector if and only $T'=\sigma
T$ for some $\sigma\in R(T)$. In this case $R(T)=R(T')$ and we say
that $T$ and $T'$ are row equivalent. Thus the vectors $e_T$ are
associated to row equivalence classes $\{T\}$ of numberings of
$\lambda$, which are called {\em tabloids} on $\lambda$. Recall that
a {\em standard tableau} on $\lambda$ is a numbering of $\lambda$
which is increasing from left to right and from top to bottom.

\begin{proposition}\label{p-SchurWeyl}
Let $\lambda$ be a partition of $N$ and let $\mathrm{dim}(V)\geq N$.
\begin{enumerate}
\item
The vectors $e_T$, where $\{T\}$ runs over tabloids on $\lambda$, form a basis
of the weight space $(V^{\otimes N})_\lambda$.
\item
The vectors $v_T=\sum_{\sigma\in C(T)}\mathrm{sign}(\sigma)
e_{\sigma T}$, where $T$ runs over the set $\mathcal T(\lambda)$ of
standard tableaux on $\lambda$, form a basis of the $S_N$-module
$W^\lambda=(V^{\otimes N})_\lambda^{\mathrm n_+}$ of primitive
vectors of weight $\lambda$.
\end{enumerate}
\end{proposition}

For proofs see, e.g., Chapter 9 of \cite{Goodman Wallach}.

A dual realization of $W^\lambda$ is also relevant. The  symmetric
bilinear form on $V$ for which the $e_i$ form an orthonormal basis
induces a symmetric non-degenerate bilinear form  $\langle\ ,\
\rangle$ on each weight space $(V^{\otimes N})_\lambda$. If
$\lambda$ is a partition of $N$, the tensors $e_T$, where $\{T\}$
runs over the set of tabloids on $\lambda$, are then an orthonormal
basis of $(V^{\otimes N})_\lambda$. Let $\tau$ be the
antiautomorphism of $\mathfrak{gl}(V)$ given by matrix transposition
with respect to the basis $e_i$. Then $\langle x\cdot
v,w\rangle=\langle v,\tau(x)\cdot w\rangle$, $x\in\mathfrak{gl}(V)$.
Moreover the bilinear form is $S_N$-invariant: $\langle\sigma\cdot
v,\sigma\cdot w\rangle=\langle v,w\rangle$, $\sigma\in S_N$, $v,w\in
V^{\otimes N}$.

\begin{proposition}\label{p-Zurbaran}
Let $\lambda$ be a partition of $N$ and let $\mathrm{dim}(V)\geq N$.
\begin{enumerate}\item[(i)]
The  form $\langle\ ,\ \rangle$
induces a non-degenerate $S_N$-invariant pairing
\[
(V^{\otimes N}/\mathfrak n_- V^{\otimes N})_\lambda\otimes
(V^{\otimes N})_\lambda^{\mathfrak n_+}\to \mathbb C.
\]
Thus we can identify  $(V^{\otimes N}/\mathfrak n_- V^{\otimes N})_\lambda$
with the dual $S_N$-module $(W^\lambda)^*$.
\item[(ii)]
The basis dual to the basis $v_T$ of $W^\lambda$ is given by
the classes of the vectors $e_T$, $T\in\mathcal T(\lambda)$ in
$(W^\lambda)^*=(V^{\otimes N}/\mathfrak n_- V^{\otimes N})_\lambda$.
\end{enumerate}
\end{proposition}

\begin{proof}
It follows from \eqref{e-hwm} and the complete reducibility of
$V^{\otimes N}$ into a direct sum of irreducible highest weight
modules that
\[
(V^{\otimes N})_\lambda=(V^{\otimes N})_\lambda^{\mathfrak n_+}
\oplus (\mathfrak n_{-}V^{\otimes N})_\lambda.
\]
Moreover this is an orthogonal direct sum with respect to the
contravariant form, since $\tau$ maps $\mathfrak n_-$ to $\mathfrak
n_+$. Therefore the pairing is well-defined and is non-degenerate.
Since $v_T$  ($T\in\mathcal T(\lambda)$) form a basis of $W^\lambda$
and $e_T$ occurs in $v_T$ with coefficient 1, we get $\langle e_T
\mod \mathfrak n_-,v_S\rangle=\delta_{T,S}$, $T,S\in \mathcal
T(\lambda)$. Thus the classes of $e_T$ form the dual basis of the
dual module.
\end{proof}

\noindent{\bf Remark.} Actually, $S_N$-modules are self-dual,
$(W^\lambda)^*\cong W^\lambda$ but the expression of the isomorphism
with respect to the bases labeled by tableaux is non-trivial. The
space of cycles in our integral formulae are more naturally
associated with $(W^\lambda)^*$.

\section{Integral representation of solutions}
\subsection{The action of the symmetric group on the solution space}
We fix a Young diagram $\lambda$ and a positive integer $m$ and keep
the notations of the introduction.

The KZ operators $D_i =\partial_i - m \sum_{j \neq i}
(s_{ij}+1)/(z_i - z_j)$ appearing in \eqref{KZ} are commuting first
order holomorphic differential operators acting on
$W^\lambda$-valued functions on the configuration space $C_N=\mathbb
C^N-\cup_{i<j}\{z_i=z_j\}$. Thus the space of holomorphic solutions
on any connected open subset $U\subset C_N$ has dimension
$\mathrm{dim}\,W^\lambda$. The symmetric group $S_N$ acts on $C_N$
and thus on the functions with values in the $S_N$-module
$W^\lambda$ by $(g\cdot\psi)(z)=g(\psi(g^{-1}\cdot z))$. The KZ
operators obey $g\cdot D_i\psi=D_{g(i)}g\cdot\psi$ for all $g\in
S_N$. Therefore $g\in S_N$ maps solutions on $U$ to solutions on
$g\cdot U$. In particular we have an action of $S_N$ on the space of
global solutions
\[
\mathit{KZ}(\lambda,m)=\{\text{holomorphic functions}\,\psi\colon
C_N\to W^\lambda \colon D_i\psi=0, \,\,i=1,\dots,N\}.
\]
In fact, all local solutions extend to global solutions:

\begin{thm}\label{t-Velasquez} \cite{O, FV} All local solutions of the
Knizhnik--Zamolodchikov equation \eqref{KZ} with $m\in\mathbb
Z_{>0}$ extend to homogeneous polynomials of degree
$m(f_2(\lambda)+(N-1)N/2)$, where $f_2(\lambda)$ is the value of the
central element $\sum_{i<j}s_{ij}$ of the group algebra of $S_N$ in
the representation $W^\lambda$. Moreover, the space of solutions
$\mathit{KZ}(\lambda,m)$ is isomorphic to $W^\lambda$ as an
$S_N$-module.
\end{thm}

The homogeneity follows directly from the equations: if $\psi$ is a
solution then $\sum z_iD_i\psi=(\sum
z_i\partial_i-m(f_2(\lambda)+(N-1)N/2))\psi=0$.

\subsection{Matsuo's integral formulae.}
The configuration spaces $C_\lambda(z)$, $z\in C_N$, of
\eqref{e-C-lambda} form a fibre bundle over $C_N$ and the action of
$S_N$ on the base lifts canonically to an action on the bundle.
Indeed, $C_\lambda(z_1,\dots,z_N)$ does not depend on the ordering
of the $z_i$. The forms $\omega_T$ of eq.~\eqref{omega} are
holomorphic differential forms on the total space that restrict to
holomorphic top differential forms $\omega_T(z)$ on each fibre
$C_\lambda(z)$. They are defined for any numberings $T$, not just
for standard tableaux and, by construction, they obey
$\omega_{gT}(g\cdot z)=\omega_T(z)$ for all $g\in S_N$, where
$gT=g\circ T$ is the natural action on the set of numberings
$T\colon \lambda\to\{1,\dots,N\}$.

Let $H_{\mathit{top}}(C_\lambda(z))_s$, $z=(z_1,\dots,z_N)\in C_N$,
be the skew-symmetric part of the homology of degree $d_\lambda$
under the action of $G_\lambda=S_{m_1}\times\cdots\times
S_{m_{n-1}}$:
\[
H_{\mathit{top}}(C_\lambda(z))_s=\{\sigma\in
H_{\mathit{top}}(C_\lambda(z))\colon g_*\sigma=(-1)^g\sigma,\quad
g\in G_\lambda\}.
\]
\begin{lemma} \label{l-omega}
If $\sigma\in H_{\mathit{top}}(C_\lambda(z))_s$ then
$\int_\sigma\omega_T(z)$ depends only on the tabloid $\{T\}$ of $T$.
\end{lemma}
Indeed, if $T$ and $T'$ differ by an element $h$ of the row group
$R(T)$, inducing a permutation $h_1$ of the set of boxes of
$\lambda$ then $\omega_T(z)=(-1)^gg^*\omega_{T'}(z)$, where $g\in
G^\lambda$ is the permutation $t_i^b\mapsto t_i^{h_1(b)}$ of the
variables (the sign comes from the volume form $dt$).

The following result can be extracted from Matsuo's paper \cite{M}.
\begin{thm}\label{t-Matsuo}(Matsuo \cite{M})
Let $\sigma\in H_{\mathit{top}}(C_\lambda(z))_s$,
$\psi_\sigma(z)=\sum_{\{T\}}\int_\sigma\omega_T(z)e_T\in (V^{\otimes
N})_\lambda$, with summation over all tabloids $\{T\}$ on $\lambda$.
\begin{enumerate}
\item[(i)] $\psi_\sigma(z)\in W^\lambda= (V^{\otimes N})^{\mathfrak
n_+}_\lambda$.
\item[(ii)] As $z$ varies in some neighbourhood of a point in $C_N$,
$\psi_\sigma(z)$ is a solution of the KZ equation \eqref{KZ}.
\end{enumerate}
\end{thm}

Recall that homology groups of neighbouring fibres of a fibre bundle
are canonically identified, so (ii) makes sense.

\begin{corollary}\label{c-psi} \
\begin{enumerate}
\item[(i)] The solution $\psi_\sigma$ of Theorem \ref{t-Matsuo} can
be written as
\[
\psi_\sigma(z)=\sum_{T\in T(\lambda)}\int_\sigma\omega_T(z)v_T.
\]
\item[(ii)]
$\psi_\sigma(g\cdot z)=g\psi_\sigma(z)$, $\sigma\in
H_{\mathit{top}}(C_\lambda(z))_s=H_{\mathit{top}}(C_\lambda(g\cdot
z))_s$.
\end{enumerate}
\end{corollary}

\begin{proof}
(i) follows from Theorem \ref{t-Matsuo}, (i) and Proposition
\ref{p-SchurWeyl}. For (ii) one uses the original expression for
$\psi_\sigma$.
\end{proof}

\subsection{Completeness}
We show here that all
solutions are given as integrals over suitable cycles. The proof is
in two parts: first we construct an $S_N$-equivariant map
$(W^\lambda)^*\to \mathit{KZ}(\lambda,m)$, defined as the integral
over the cycles $\sigma_T$. Since $(W^\lambda)^*$ is irreducible and
of the right dimension, it then suffices to check that the map is
non-zero, which can be done by an asymptotic analysis.

We start by describing the action of the symmetric group on cycles.

\begin{lemma}\label{l-sigma}
For any numbering $T$ of $\lambda$ and $z=(z_1,\dots,z_N)\in C_N$,
let $\sigma_T=\sigma_T(z)\in H_{\mathit{top}}(C_\lambda(z))_s$ be
the homology class defined in the introduction as the
skew-symmetrization of the image of the fundamental class by a map
$(S^1)^{d_\lambda}\to C_\lambda(z)$.
\begin{enumerate}
\item[(i)] For all $g\in S_N$, $\sigma_{g T}(g\cdot z)=\sigma_{T}(z)$.
\item[(ii)] If $T$ and $T'$ are numberings of $\lambda$ differing by
a row permutation, then $\sigma_T(z)=\sigma_{T'}(z)$. Thus
$\sigma_T(z)$ depends only on the tabloid of $T$.
\end{enumerate}
\end{lemma}

\begin{proof} (i) holds by construction. The proof of (ii)
is the same as the proof of Lemma \ref{l-omega}. This time the sign
comes from the change of orientation.
\end{proof}

Thus we get a map
\begin{eqnarray*}
\Psi^\lambda_m\colon (V^{\otimes
N})_\lambda&\to&\mathit{KZ}(\lambda,m)
\\
e_T&\mapsto&\int_{\sigma_T}\omega=\sum_{T'\in
T(\lambda)}\psi_{T,T'}v_{T'}.
\end{eqnarray*}
It is well-defined since $\sigma_T$, just as $e_T$, depends only on
the tabloid of $T$.

\begin{lemma}\label{p-Psi-m}
 The map $\Psi^\lambda_m$ is $S_N$-equivariant.
\end{lemma}

This  follows from Lemma \ref{l-sigma} and Corollary \ref{c-psi}.

To prove that the map is non-zero we will use the following key technical lemma. Let $(r(b),c(b))$ be the coordinates (row and column number) of the box $b$ in the Young diagram $\lambda$.
\begin{lemma}\label{l-nonzero}
Let $T$ be the standard tableau mapping the $k$th box, counted from
left to right and top to bottom, to $k$. Then $\psi_{T,T}(z)$ is not
identically zero. The leading term for $0\ll |z_1|\ll
|z_2|\ll\cdots\ll |z_N|$ is
\[
\psi_{T,T}(z)\sim C\, \prod_{b\in\lambda}
z_{T(b)}^{m(T(b)-1+c(b)-r(b))}+\cdots,
\]
for some $C\neq 0$.
\end{lemma}

\begin{proof}
We show that, as $z_N\to\infty$,
\[\psi_{T,T}(z_1,\dots,z_N)=C'
z_N^{m(N-1+\lambda_n-n)}(\psi_{T',T'}(z_1,\dots,z_{N-1})+O(z_N)),
\]
where $T'$ is the standard tableau with $N-1$ boxes obtained from
$T$ by removing the last box and $C'$ is a non-zero combinatorial
factor. Since $\psi_{\bx{1},\bx{1}}=1$ for the tableau with one box,
this gives an inductive proof of the claim.

Let $\lambda'$, the shape of $T'$, be $\lambda$ without the last
box. Then
\begin{eqnarray*}
 \Phi_{\lambda}&=&\Phi_{\lambda'}\prod_{k=1}^{N-1}(z_k-z_N)^2(t_1^{\bx{k}}-z_N)^{-1}
 \prod_{s=1}^{n-1}(t_s^{\bx{N}}-t_{s-1}^{\bx{N}})^{-1}
 \\
 &&\prod_{s=1}^{n-1}\left(
 \prod_{\stackrel{k<N}{r(\bx{k})>s}}(t_s^{\bx{k}}-t_s^{\bx{N}})^2
 \prod_{\stackrel{k<N}{r(\bx{k})>s+1}}(t_{s+1}^{\bx{k}}-t_{s}^{\bx{N}})^{-1}
 \prod_{\stackrel{k<N}{r(\bx{k})>s-1}}(t_s^{\bx{N}}-t_{s-1}^{\bx{k}})^{-1}\right),
\end{eqnarray*}
where $\bx k=T^{-1}(k)$ is the box of $\lambda$ labeled by $k$ and
we set $t_0^{\bx{k}}=z_k$. Also,
\[
\phi_{T}=\phi_{T'}\prod_{s=1}^{n-1}(t_s^{\bx{N}}-t_{s-1}^{\bx{N}})^{-1}.
\]
The leading term as $z_N\to\infty$ in $\psi_{T,T}$ is obtained when
the variables $t^{\bx{N}}_s$ run over circles around $z_N$. With the
variable substitution $t^{\bx{N}}_s=z_N+\tau_1+\cdots+\tau_s$, the
leading term as $z\to\infty$ of the integration of $\omega_T$ over
the variables $t^{\bx{N}}_s$ is
\[
\pm\omega_{T'}z_N^{m(N-1)}\mathrm{res}_{\tau_{n-1}=0}\cdots
\mathrm{res}_{\tau_1=0}\Omega,
\]
where
\begin{gather*}
\Omega=
 \prod_{s=1}^{n-2}(z_N+\tau_1+\cdots+\tau_s)^{m(2m_s-m_{s-1}-m_{s+1})}
 \tau_{s}^{-m-1}d\tau_s
 \\
 \cdot(z_N+\tau_1+\cdots+\tau_{n-1})^{m(2m_{n-1}-m_{n-2}-1)}
 \tau_{n-1}^{-m-1}d\tau_{n-1}
 \\
=
\prod_{s=1}^{n-2}(z_N+\tau_1+\cdots+\tau_s)^{m(\lambda_{s+1}-
\lambda_s)}\tau_{s}^{-m-1}d\tau_s
 \\
 \cdot(z_N+\tau_1+\cdots+\tau_{n-1})^{m(\lambda_{n}-\lambda_{n-1}-1)}
 \tau_{n-1}^{-m-1}d\tau_{n-1}.
\end{gather*}
The residues can be computed explicitly. Such expressions give a
non-zero result if the total power of all factors containing any
given $\tau_s$ is {\em negative}. This power is
$(\lambda_n-\lambda_s-1)m$ which is indeed negative for all
$s=1,\dots,n-1$.
\end{proof}

\begin{thm}\label{t-1}
The map $\Psi^\lambda_m\colon (V^{\otimes N})_\lambda\to
\mathit{KZ}(\lambda,m)$ is an epimorphism of $S_N$-modules with
kernel $(\mathfrak n_-V^{\otimes N})_\lambda$ and therefore defines
an isomorphism $(W^\lambda)^*\to \mathit{KZ}(\lambda,m)$. In
particular, the images $\psi_T$ of the basis vectors $[e_T]$,
$T\in\mathcal T(\lambda)$, of $(W^\lambda)^*=(V^{\otimes
N}/\mathfrak n_-V^{\otimes N})_\lambda$, form a basis of the space
of solutions of the KZ equation.
\end{thm}

\begin{proof} By Lemma \ref{p-Psi-m} and Lemma \ref{l-nonzero},
$\Psi^\lambda_m$ is a non-zero homomorphism from the $S_N$-module
$(V^{\otimes N})_\lambda$ to the $S_N$-module
$\mathit{KZ}(\lambda,m)$. By Theorem \ref{t-Velasquez}
$\mathit{KZ}(\lambda,m)$ is isomorphic to the irreducible
$S_N$-module $W^\lambda$. Since the image of a homomorphism is a
submodule, it follows that the map is surjective. On the other hand,
by \eqref{e-Ribera}, $W^\lambda\simeq (W^\lambda)^*$ occurs in
$(V^{\otimes N})_\lambda$ with multiplicity
$\mathrm{dim}\,M^\lambda_\lambda=1$ and the claim follows from
Proposition \ref{p-Zurbaran}.
\end{proof}

Another interesting corollary of Lemma \ref{l-nonzero} is the
following classical formula.
\begin{proposition} (Frobenius \cite{F})
Let $f_2(\lambda)$ be the value of the central element
$\sum_{i<j}s_{ij}$ in the representation $W^\lambda$. Then
\begin{equation}
\label{Frob} f_2(\lambda)=\sum_{b\in\lambda}(c(b)-r(b)),
\end{equation}
where as before $r(b)$ and $c(b))$ are respectively the row and
column coordinates of the box $b$ in the Young diagram $\lambda$.
\end{proposition}
To prove this we recall that the degree of the polynomial solutions
from $\mathit{KZ}(W,m)$ is equal to the value of the central element
$\sum_{i<j} m (s_{ij}+1)$ in the irreducible representation $W$ (see
\cite{FV}). Comparing this with the leading term of the solution
from Lemma \ref{l-nonzero} and taking into account that $\sum_{i<j}
1 = \sum_{b\in\lambda}( T(b)-1) = \frac{N(N-1)}{2}$ we come to
Frobenius formula (\ref{Frob}).

\subsection{Integrality}
It is well-known (and clear from Proposition \ref{p-SchurWeyl}) that
$(v_T)_{T\in\mathcal T(\lambda)}$ is an integral basis of
$W^\lambda$, i.e., $W_{\mathbb
Z}^\lambda=\oplus_{T\in\mathcal{T}(\lambda)} \mathbb Z\,v_T$ is a
module over the group ring $\mathbb Z S_N$.

\begin{thm}\label{t-2}
The functions $\psi_{T,T'}$ are homogeneous polynomials in
$z_1,\dots,z_N$ with integer coefficients. Thus
\[
\Psi^\lambda_m(z_1,\dots,z_N)=\sum_{T,T'\in\mathcal
T(\lambda)}\psi_{T,T'}(z_1,\dots,z_N)\,v_T\otimes v_{T'}\in \mathbb
W_{\mathbb Z}^\lambda\otimes_\mathbb Z W_{\mathbb
Z}^\lambda[z_1,\dots,z_N].
\]
Moreover $\Psi^\lambda_m$ is $S_N$-invariant: $\Psi^\lambda_m(g\cdot
z)=g\otimes g\,\Psi^\lambda_m(z)$, for all $g\in S_N$.
\end{thm}

The invariance is a rephrasing of the homomorphism property of
$\Psi_m^\lambda$ of Theorem \ref{t-1}. The integrality follows from
repeated application of the following elementary

\begin{lemma}
Let $m_i$ be any integers and $w_1,\dots,w_k$ distinct complex numbers.
Then, for any contour $\gamma$ in the complex plane not passing
through $w_1,\dots,w_k$,
\[
\frac1{2\pi i}\oint_\gamma\prod_{i=1}^k(t-w_i)^{m_i}dt =
\sum_{i<j}c_{ij}(w_i-w_j)^{\ell_{ij}},
\]
for some integers $c_{ij}$, $\ell_{ij}$.
\end{lemma}

Theorems \ref{t-1} and \ref{t-2} imply the statements of Theorem
\ref{t-main}.

\noindent{\bf Example.} We demonstrate here our formulae in the simplest non-trivial example when $N=3$ and $\lambda= (2,1),$ which corresponds to the usual two-dimensional representation of $S_3$. In this case there are two standard tableaux:

\centerline { $T=$
\begin{picture}(30,30)(0,50)
\put( 0,55){\framebox(14.50,14.60)[c]{1}}
\put(15,55){\framebox(14.60,14.60)[c]{2}}
\put(0,40){\framebox(14.60,14.60)[c]{3}}
\end{picture},
\qquad $S=$
\begin{picture}(30,30)(0,50)
\put( 0,55){\framebox(14.50,14.60)[c]{1}}
\put(15,55){\framebox(14.60,14.60)[c]{3}}
\put(0,40){\framebox(14.60,14.60)[c]{2}}
\end{picture}.
}

\bigskip

\noindent The corresponding primitive vectors are
$v_T=\epsilon_3-\epsilon_1$, $v_S=\epsilon_2-\epsilon_1$.

The residues can be computed explicitly and one obtains a basis of
$\mathit{KZ}(\lambda,m)$:
\begin{eqnarray*}
\psi_1(z_1,z_2,z_3)&=&z_{23}^{2m}\sum_{k=0}^m d_{m,k}
\left((m-k)v_T+kv_S\right)z_{12}^{m-k}z_{13}^k,
\\
\psi_2(z_1,z_2,z_3)&=&z_{13}^{2m}\sum_{k=0}^m(-1)^{m-k}
d_{m,k}\left((m-k)v_T-mv_S\right)z_{12}^{m-k}z_{23}^k,
\end{eqnarray*}
where we abbreviate $z_i-z_j=z_{ij}$ and
\[
d_{m,k}=-\frac1m{-m\choose k}{-m\choose m-k}.
\]

\section{Duality $m \leftrightarrow -m$ and intersection pairing}

To apply our results to the case of negative integers $m$ we can use the isomorphism between the spaces of solutions
$\mathit{KZ}(V,m)$ and $\mathit{KZ}(V\otimes\mathrm{Alt},-m)$ mentioned in the Introduction.

\begin{lemma}
For any $\psi\in\mathit{KZ}(V,m)$ the function
$$\phi=\prod_{i<j}(z_i-z_j)^{-2m}
\psi\otimes\epsilon$$ belongs to the space of solutions $\mathit{KZ}(V\otimes\mathrm{Alt},-m)$. \end{lemma}

Indeed from the KZ equation (\ref{KZ}) we have
$$
\partial_i \phi =
(-2m) (\sum_{j \neq i}^N \frac{1}{z_i - z_j})
\prod_{i<j}(z_i-z_j)^{-2m}\psi\otimes\epsilon + m
\prod_{i<j}(z_i-z_j)^{-2m}(\sum_{j \neq i}^N \frac{s_{ij}+1}{z_i -
z_j}\psi)\otimes\epsilon = $$ $$ (-m)
\prod_{i<j}(z_i-z_j)^{-2m}\sum_{j \neq i}^N \frac{s_{ij}+1}{z_i -
z_j}(\psi\otimes\epsilon) = (-m) \sum_{j \neq i}^N
\frac{s_{ij}+1}{z_i - z_j}\phi, $$ since
$s_{ij}(\phi\otimes\epsilon) = - (s_{ij}\phi)\otimes\epsilon.$

{\bf Remark.} The pre-factor $\prod_{i<j}(z_i-z_j)^{-2m}$
will disappear if we consider the KZ equation in the form
$$\partial_i \psi = m \sum_{j \neq i}^N \frac{s_{ij}}{z_i - z_j}
\psi, \quad i=1,\dots,N,$$ which in general has rational
solutions. The duality between $m$ and $-m$ for similar systems
was used in \cite{FV-1} to explain the shift operator for the
Calogero-Moser systems.

It is well-known that the involution $V \rightarrow
V\otimes\mathrm{Alt}$ corresponds to the {\it transposition} of
the Young diagram $\lambda \rightarrow \lambda'.$ Thus we have
established an isomorphism
\begin{equation}
\label{isoi} \mathit{KZ}(\lambda,m) \approx
\mathit{KZ}(\lambda',-m).
\end{equation}
Note that the configuration spaces
$C_{\lambda} (z_1, \dots, z_N)$ and $C_{\lambda'} (z_1, \dots,
z_N)$  in general are quite different (in particular, they have
different dimensions), so the structure of the integral formulas
for the solution of the KZ equations for a given Young diagram,
which we get in this way, substantially depends on the sign of
$m.$

It turns out that there is a link between the spaces of KZ solutions with the {\it same} Young diagram:
\begin{equation}
\label{dua1}
j: \mathit{KZ}(\lambda,m) \approx \mathit{KZ}(\lambda,-m)^*.
\end{equation}
More precisely, there exists a natural pairing
\begin{equation}
\label{dua2} \mathit{KZ}(V, m)\times \mathit{KZ}(V^*,
-m)\to\mathbb C,
\end{equation}
where $V^*$ is the dual space to $V$ with the natural action of $S_N.$
It is defined by the following lemma.
\begin{lemma}\label{pairing}
Let $\langle\ ,\ \rangle$ denote the canonical pairing between $V$
and $V^*.$ Then for any two solutions $\psi(z_1,\dots,z_N) \in
\mathit{KZ}(V, m)$ and $\phi(z_1,\dots,z_N) \in \mathit{KZ}(V^*,
-m)$ the product
$\langle\psi(z_1,\dots,z_N),\phi(z_1,\dots,z_N)\rangle$ is
independent of $z_1, \dots,z_N$ and thus defines a non-degenerate
canonical pairing (\ref{dua2}).
\end{lemma}

The proof is a straightforward check using the KZ equations
(\ref{KZ}):
\[
\partial_i \langle\psi(z_1,\dots,z_N),\phi(z_1,\dots,z_N)\rangle =  \langle m\sum_{j \neq i}^N \frac{s_{ij}+1}{z_i - z_j}
\psi, \phi \rangle +  \langle \psi, (-m)\sum_{j \neq i}^N \frac{s_{ij}+1}{z_i - z_j}\phi \rangle = 0,
\]
since $\langle s_{ij} v, w \rangle = \langle v, s_{ij} w \rangle$ for any $v \in V, w \in V^*.$
The same calculation holds of course for any Coxeter group.

Note that the $S_N$-module $V^*$ is isomorphic to $V$ and in the
irreducible case the isomorphism is {\it almost canonical} in the
sense that it is unique up to a scaling factor. This leads to an
isomorphism (\ref{dua1}) and allows us to find the solutions from
$\mathit{KZ}(\lambda,-m)$ as follows.

Let us choose any basis $e_1, \dots, e_M$ in $V= W^{\lambda}$ and a
basis of solutions $\psi_{\alpha} = \sum_{i=1} ^M \psi_{\alpha}^i
(z_1, \dots, z_N) e_i, $ in $\mathit{KZ}(\lambda,m).$ Let
$\Phi_{\lambda,m}(z_1,\dots,z_N) = \| \psi_{\alpha}^i (z_1, \dots,
z_N)\|$ be the corresponding $M \times M$ {\it fundamental matrix}
of $\mathit{KZ}(\lambda,m)$. Let $e^1, \dots, e^M$ be the dual basis
in $V^*: \langle e^i , e_j\rangle=\delta^i_j.$ We are looking now
for a fundamental matrix
$\Phi_{\lambda,-m}(z_1,\dots,z_N)=\|\phi_\beta^j(z_1,\dots,z_N)\|$
for $\mathit{KZ}(\lambda,-m)$, given by a basis of solutions
$\phi^\beta=\sum_{j=1}^M\phi^\beta_j(z_1,\dots,z_N)e^j$. 
From Lemma \ref{pairing} it follows that one can choose the basis of solutions in such a way that
\begin{equation}\label{invt}
 \sum_{i=1}^M\psi_\alpha^i(z_1,\dots,z_N)\phi^\beta_i(z_1,\dots,z_N)
=\delta_\alpha^{\beta}.
\end{equation}

\begin{proposition}
A fundamental matrix for $\mathit{KZ}(\lambda,-m)$ can be found as
the tranposed inverse matrix to the fundamental matrix of
$\mathit{KZ}(\lambda,m):$
$$
\Phi_{\lambda,-m}(z_1,\dots,z_N) =
(\Phi_{\lambda,m}(z_1,\dots,z_N)^{-1})^{T}.
$$
The determinant of the fundamental matrix is given by
\begin{equation}\label{det}
\det\Phi_{\lambda,m}(z_1,\dots,z_N)= C \prod_{i<j}(z_i-z_j)^{2m
d_+(\lambda)},
\end{equation} where $C=C(\lambda,m)$ is a non-zero constant and $d_+(\lambda)= \dim W^{\lambda}_+$ is the dimension of
the fixed subspace of reflection $s_{ij}$ acting in the
representation $W^{\lambda}.$
\end{proposition}

The first part is equivalent to (\ref{invt}).
To prove the formula for the determinant we can use the standard fact that if
matrix $\Phi$ satisfies the differential equation $\Phi' = A \Phi$
then its determinant satisfies the equation $\det \Phi ' = tr A
\det\Phi.$ Applying this to the KZ equation (\ref{KZ}) and using
the fact that $tr(s_{ij}+1) = 2d_+(\lambda)$ we have the result.
The formula (\ref{det}) shows that the singularities of
$\Phi_{\lambda,-m}(z_1,\dots,z_N)$ are located on the hyperplanes
$z_i=z_j$, which of course follows also from the previous
considerations.

 In the rest of this section we discuss the topological interpretation of the duality (\ref{dua2}) as
intersection pairing between certain homology groups. It is based
on the integral formula for the solutions of the KZ equation found
in our previous work \cite{FV} (see section 4.5).

We restrict ourselves with the special case of the reflection
representation of $S_N$, which is the standard $(N-1)$-dimensional
irreducible representation on the hyperplane $x_1+\dots+x_n=0$ in
$\mathbb C^N$ defined by permutation of coordinates. This
representation is isomorphic to $W^\lambda$ with
$\lambda=(N-1,1)$. For positive $m$ our integrals giving the
solutions are one-dimensional and, in terms of the standard basis
$\epsilon_b$ of $\mathbb C^N$, they have the form
\[
\psi_a=\prod_{1\leq i<j\leq
N}(z_i-z_j)^{2m}\,\mathrm{res}_{t=z_a}\prod_{i=1}^N(t-z_i)^{-m}\sum_{b=1}^N\frac1{t-z_b}\epsilon_b\,dt,
\quad a=1,\dots,N.
\]
They obey the relation $\psi_1+\cdots+\psi_N=0$ and
$\psi_1,\dots,\psi_{N-1}$ form a basis of the solution space
$\mathit{KZ}(\lambda,m)$.

A different integral representation for the solution space
$\mathit{KZ}(\lambda,-m)$ for positive $m$ was found in \cite{FV},
where it was shown that
\[
\phi_a=\prod_{1\leq i<j\leq N}(z_i-z_j)^{-2m}\,\int_{z_a}^{z_{N}}
\prod_{i=1}^N(t-z_i)^{m}\sum_{b=1}^N\frac1{t-z_b}\epsilon_b\,dt,\qquad
a=1,\dots,N-1
\]
give a basis in $\mathit{KZ}(\lambda,-m).$ In particular, in the
case $N=3$ we have after explicit evaluation of the integrals the
following basis:
\begin{eqnarray*}
\phi_1(z_1,z_2,z_3)&=&z_{23}^{-2m}\sum_{k=0}^m(-1)^{m+k} d'_{m,k}
\left((-m-k) v_T+kv_S\right)z_{12}^{-m-k}z_{13}^k,
\\
\phi_2(z_1,z_2,z_3)&=&z_{13}^{-2m}\sum_{k=0}^m d'_{m,k}
\left((-m-k) v_T+mv_S\right)z_{12}^{-m-k}z_{23}^k,
\end{eqnarray*}
where
\[
d'_{m,k}={m\choose k}\frac{(m-1)!(m+k-1)!}{(2m+k)!}.
\]
\medskip

Thus, in a more invariant geometric terms, for $\lambda=(N-1,1)$ and $m\in\mathbb Z_{>0}$, we
have two maps
\[
\psi\colon
H_1(\mathbb C\smallsetminus\{z_1,\dots,z_N\})\to
W^\lambda,\qquad \phi\colon H_1(\mathbb C,\{z_1,\dots,z_N\})\to
W^\lambda,
\]
sending horizontal sections for the Gauss--Manin connection to
solutions in $\mathit{KZ}(\lambda,m)$ and in
$\mathit{KZ}(\lambda,-m)$, respectively. The map $\phi$ induces an
isomorphism between the complexification of the space of horizontal
relative cycles and $\mathit{KZ}(\lambda,-m)$, whereas $\psi$ has a
one-dimensional kernel spanned by a cycle surrounding all the points
$z_i$. This kernel is exactly the complexification of the left
kernel of the intersection pairing
\[
H_1(\mathbb C\smallsetminus\{z_1,\dots,z_N\})\times H_1(\mathbb
C,\{z_1,\dots,z_N\})\to \mathbb Z,
\]
and the right kernel is trivial. Since the intersection pairing of
cycles is preserved by the Gauss--Manin connection we obtain a
non-degenerate $S_N$-invariant pairing $\mathit{KZ}(\lambda,m)\times
\mathit{KZ}(\lambda,-m)\to\mathbb C$.

The claim is that this pairing is proportional to the one described in Lemma \ref{pairing}.

\begin{proposition} Let $\lambda=(N-1,1)$ and $m>0$.
The pairing (\ref{dua2}) of solution spaces $\mathit{KZ}(\lambda,m)$
and $\mathit{KZ}(\lambda,-m)$ is proportional to the image of the
intersection pairing $(\ \cdot\ )$. More precisely, let $\sigma\in
H_1(\mathbb C\smallsetminus\{z_1,\dots,z_N\}),\quad \tau\in
H_1(\mathbb C,\{z_1,\dots,z_N\})$ vary with the points $z_i$ as
horizontal sections, then
\begin{gather*}
\langle\psi_\sigma(z_1,\dots,z_N),\phi_\tau(z_1,\dots,z_N)\rangle=
C_{N}\frac1m\,(\sigma\cdot\tau),
\\
\sigma\in H_1(\mathbb
C\smallsetminus\{z_1,\dots,z_N\}),\quad \tau\in H_1(\mathbb
C,\{z_1,\dots,z_N\}),
\end{gather*}
for some constant $C_{N}\neq 0$ depending on the normalization of
the isomorphism $(W^\lambda)^*\to W^\lambda$.
\end{proposition}

\begin{proof}
The proof follows from Schur's lemma, except for the
determination of the constant of proportionality of the two
pairings. To compute it, we consider two special cycles, namely a
small circle around $z_1$ and a path from $z_1$ to $z_N$. These
cycles have intersection number $-1$ and the corresponding solutions
are $\psi_1$ and $\phi_1$, respectively. It is sufficient to compute
the pairing when $z_1=0$ in the limit $z_N\to 0$, where also
$\phi_1$ is regular. After the change of variable $t=z_N\tau$ we
obtain
\begin{gather*}
\psi_1(0,z_2,\dots,z_{N-1},0)=F_m\,\mathrm{res}_{\tau=0}\tau^{-m}(\tau-1)^{-m}\left(\frac{\epsilon_1}\tau+
\frac{\epsilon_N}{\tau-1}\right)d\tau,
\\
\phi_1(0,z_2,\dots,z_{N-1},0)=F^{-1}_m\int_0^1\tau^{m}(\tau-1)^{m}\left(\frac{\epsilon_1}\tau+
\frac{\epsilon_N}{\tau-1}\right)d\tau,
\end{gather*}
for some function $F_m$ of $z_2,\dots,z_{N-1}$. By integrating by
parts we see that the coefficient of $\epsilon_N$ is minus the
coefficient of $\epsilon_1$ (as it should be since the solutions
take values in primitive vectors). The result of the calculation is
\begin{gather*}
\psi_1(0,z_2,\dots,z_{N-1},0)=
F_m\,(-1)^m\frac{(2m-1)!}{m!(m-1)!}(\epsilon_1-\epsilon_N),
\\
\phi_1(0,z_2,\dots,z_{N-1},0)=
F_m^{-1}(-1)^m\frac{m!(m-1)!}{(2m)!}(\epsilon_1-\epsilon_N).
\end{gather*}
If we normalize the pairing $W^\lambda\times W^\lambda\to \mathbb C$
defining the isomorphism between $W^\lambda$ and $(W^\lambda)^*$ so
that the basis $\epsilon_a$ is orthonormal, we obtain
$\langle\psi_1,\phi_1\rangle=1/m$.
\end{proof}

Note that since the cycles defining the bases $\psi_a$ and $\phi_a$
are dual (up to sign) with respect to the intersection pairing, we
deduce
\[
\langle\psi_a(z_1,\dots,z_N),\phi_b(z_1,\dots,z_N)\rangle
=-C_{N}\frac1m\,\delta_{a,b},
\]
where $C_N=-1$ with the choice of normalization described in the
proof.

\medskip

\noindent {\bf Remark.}
We would like to mention that an extension of these results to the
case of general representations should involve a suitable
replacement of relative homology. The singularity as integration
variables approach each other causes difficulties with the naive
generalization.

\section*{Acknowledgements} We are grateful to A. Okounkov, who
helped us to trace back the Frobenius formula, to O. Chalykh, who attracted our attention to Opdam's lecture notes \cite{O}, 
and to A. Chervov for stimulating discussions.

 This work has been partially supported by the
European Union through the FP6 Marie Curie RTN ENIGMA (Contract
number MRTN-CT-2004-5652). It was made possible by a visit of A.
V. to FIM at ETH, Zurich with the support of the MISGAM programme
of the European Science Foundation, and a visit of the authors to
the Newton Institute for Mathematical Sciences, Cambridge. We
thank these institutions for their support and hospitality.


\begin{thebibliography}{8}

\bibitem{FeiV}
M. Feigin and A. P. Veselov. {\it Quasi-invariants of Coxeter groups
and m-harmonic polynomials.} Int. Math. Res. Notices (2002), no. 10,
521-545.

\bibitem{FV}
G. Felder and A. P. Veselov.  {\it Action of Coxeter groups on
m-harmonic polynomials and KZ equations.} Moscow Math Journal, 2003,
vol.3, no.4, 1269-1291.

\bibitem{FV-1}
G. Felder and A. P. Veselov.  {\it Shift operators for the quantum
Calogero-Sutherland problems via Knizhnik-Zamolodchikov equation.}
Commun. Math. Phys. 160 (1994), 259-273.

\bibitem{F}
F. G. Frobenius. {\it \"Uber die Charaktere der symmetrischen
Gruppe.} Sitz. K\"onig. Preuss. Akad. Wissenschaften zu Berlin
(1900), 516-534.

\bibitem{Fulton}
W. Fulton. {\it Young Tableaux, with Applications to Representation
Theory and Geometry.} Cambridge University Press, 1997

\bibitem{Goodman Wallach}
R. Goodman and N. R. Wallach. {\it Representations and Invariants of
the Classical Groups.} Cambridge University Press, 1998

\bibitem{KK}
A. Kazarnovski-Krol. {\it Cycles for asymptotic solutions and the
Weyl group.} The Gelfand Mathematical Seminars, 1993--1995,
123--150, Gelfand Math. Sem., Birkh\"{a}user Boston, Boston, MA,
1996.

\bibitem{Mac}
I. G. Macdonald. {\it Symmetric functions and Hall polynomials.}
Second edition, Oxford University Press 1999.

\bibitem{M}
A. Matsuo. {\it An application of Aomoto-Gelfand hypergeometric
functions to the SU(n) Knizhnik-Zamolodchikov equation.} Commun.
Math. Phys. 134 (1990), 65-77.

\bibitem{O}
E. Opdam. {\it Complex reflection groups and fake degrees.}
arXiv: math.RT/9808026.

\bibitem{Varchenko}
A. Varchenko. {\it Hypergeometric Functions and Representation
Theory of Lie Algebras and Quantum Groups.} Advanced Series in
Mathematical Physics, Vol.21, World Scientific (1995)
\end{thebibliography}
\end{document}